\documentclass{gtmon_a}
\pdfoutput=1

\proceedingstitle{Groups, homotopy and configuration spaces (Tokyo
  2005)}
\conferencestart{5 July 2005}
\conferenceend{11 July 2005}
\conferencename{Groups, homotopy and configuration spaces, 
                in honour of Fred Cohen's 60th birthday}
\conferencelocation{University of Tokyo, Japan}

\editor{Norio Iwase}
\givenname{Norio}
\surname{Iwase}

\editor{Toshitake Kohno}
\givenname{Toshitake}
\surname{Kohno}

\editor{Ran Levi}
\givenname{Ran}
\surname{Levi}

\editor{Dai Tamaki}
\givenname{Dai}
\surname{Tamaki}

\editor{Jie Wu}
\givenname{Jie}
\surname{Wu}

\title[On the L-S category of symmetric spaces of classical type]
      {On the Lusternik--Schnirelmann category of\\
       symmetric spaces of classical type}

\author[M Mimura]{Mamoru Mimura}
\givenname{Mamoru}
\surname{Mimura}
\address{Department of Mathematics\\
         Faculty of Science\\
         Okayama University\\\newline
         Okayama 700-8530, Japan}
\email{mimura@math.okayama-u.ac.jp}
\urladdr{}

\author[K Sugata]{Kei Sugata}
\givenname{Kei}
\surname{Sugata}
\email{sugata@math.okayama-u.ac.jp}
\urladdr{}


\arxivreference{}

\volumenumber{13}
\issuenumber{}
\publicationyear{2008}
\papernumber{15}
\startpage{323}
\endpage{334}
\doi{}
\MR{}
\Zbl{}

\received{18 February 2006}
\revised{4 September 2006}
\accepted{11 October 2006}
\published{25 February 2008}
\publishedonline{25 February 2008}
\proposed{}
\seconded{}
\corresponding{}
\version{}

\makeatletter
\def\cnewtheorem#1[#2]#3{\newtheorem{#1}{#3}[section]
\expandafter\let\csname c@#1\endcsname\c@thm}

\AtBeginDocument{\let\bar\wbar\let\tilde\wtilde\let\hat\what}

\allowdisplaybreaks

\newtheorem*{Theorem}{Theorem}
\newtheorem{thm}{Theorem}[section]
\cnewtheorem{lem}[thm]{Lemma}
\cnewtheorem{prop}[thm]{Proposition}
\numberwithin{equation}{section}

\newtheorem{thmm}{Theorem}

\makeatother  

\renewcommand{\u}{\mathfrak{u}}
\AtBeginDocument{\renewcommand{\t}[1]{{}^t\!{#1}}}

\DeclareMathOperator{\cat}{cat}
\DeclareMathOperator{\cupl}{cup}
\DeclareMathOperator{\tr}{tr}

\begin{document}

\begin{htmlabstract}
We determine the Lusternik&ndash;Schnirelmann category of the irreducible,
symmetric Riemann spaces SU(n)/SO(n) and SU(2n)/Sp(n) of type AI and AII
respectively.
\end{htmlabstract}

\begin{abstract}    
We determine the Lusternik--Schnirelmann category of the irreducible,
 symmetric Riemann spaces $SU(n)/SO(n)$ and $SU(2n)/Sp(n)$ of type AI and AII
 respectively.
\end{abstract}
\begin{asciiabstract}
We determine the Lusternik-Schnirelmann category of the irreducible,
 symmetric Riemann spaces SU(n)/SO(n) and SU(2n)/Sp(n) of type AI
 and AII respectively.
\end{asciiabstract}

\maketitle


\section{Introduction}

For a topological space $X$, the Lusternik--Schnirelmann category, L-S category
 for short and denoted by $\cat(X)$, is defined to be the least integer $n$
 such that there exists an open covering $\{A_1,\dots,A_{n+1}\}$ of $X$ with
 each $A_i$ contractible in $X$.
This homotopy invariant is known to be related to various problems; for
 instance, some geometric applications can be found in Korba{\v{s}} and 
Sz{\H{u}}cs \cite{K1}.

First of all we recall a theorem due to Singhof \cite{S}:

\begin{Theorem}$
\cat(SU(n)) = n-1
$
\end{Theorem}

The purpose of this note is to prove the following theorem along the line of
 idea of the proof of Singhof's theorem. 

\begin{thm}\label{main}
\ \\
\ $(1)$\qua $\cat(SU(n)/SO(n)) = n-1$\\
\ $(2)$\qua $\cat(SU(2n)/Sp(n)) = n-1$
\end{thm}

One can prove the following theorem by the entirely similar method.

\begin{thmm}
\ \\
\ $(1)$\qua $\cat(U(n)/O(n)) = n$\\
\ $(2)$\qua $\cat(U(2n)/Sp(n)) = n$
\end{thmm}

Observe that (1) of \fullref{main} for $n=4$ improves the estimate of the
 L-S category of the oriented Grassmann manifold
 $\tilde{G}_{6,3} = SO(6)/(SO(3) \times SO(3))$ given by
 Korba{\v{s}} \cite[Corollary C (a)]{K2}.

Let
$
J = \left(
  \begin{array}{@{}c@{\,}c@{}}
       O & -E_n \\
     E_n & O
  \end{array}
\right) \in SU(2n)
$,
 where $E_n$ denotes the $n \times n$ identity matrix.

We need the following lemma to give a proof of our result.

\begin{lem}\label{homeo}
There are matrix representations\textup{:}\\
\ $(1)$\qua $SU(n)/SO(n) = \{ X \in SU(n) \mid \t{X} = X \}$\\
\ $(2)$\qua $SU(2n)/Sp(n) = \{ X \in SU(2n) \mid \t{X} = JX\t{J} \}$
\end{lem}

By \fullref{homeo}, we can regard $SU(n)/SO(n)$ and $SU(2n)/Sp(n)$ as
 subspaces of $SU(n)$ and $SU(2n)$ respectively.

The paper is organized as follows.
In Section 2 we will prove (1) of \fullref{main}.
In Section 3 we will prove (2) of \fullref{main}.
In Section 4 we study the L-S category of the irreducible symmetric
 Riemann spaces of classical type other than AI and AII.
We will give a proof of \fullref{homeo}, which may be a folklore,
 in the Appendix just for completeness.

We thank J~Korba\v{s} for reading the manuscript, M~Yasuo for giving
 information on these homogenous spaces and also N~Iwase and T~Nishimoto for
 useful conversations.

\section{L-S category of $SU(n)/SO(n)$}\label{SU(n)/SO(n)}

In this section, we will prove (1) of \fullref{main}.
The mod 2 cohomology ring of $SU(n)/SO(n)$ is given as follows (see for
 example Mimura and Toda \cite{MT}):
\[
H^*(SU(n)/SO(n);\Z/2) = \Lambda(x_2,x_3,\dots,x_n),
\]
 where $\Lambda$ denotes exterior algebra.
Since the cup-length gives a lower bound of the L-S category (see for example
 Whitehead \cite{Wh}), we have
\[
n-1 = \cupl_{\Z/2}(SU(n)/SO(n)) \leq \cat(SU(n)/SO(n)).
\]
Thus in order to determine $\cat(SU(n)/SO(n))$, it is sufficient to show the
 following proposition.

\begin{prop}
$
\cat(SU(n)/SO(n)) \leq n-1
$
\end{prop}
\proof
Let $\lambda_1,\lambda_2,\dots,\lambda_n$ be different complex numbers with
 $|\lambda_r|=1$ such that $\lambda_1\lambda_2\cdots\lambda_n \neq 1$.
For $1 \leq r \leq n$ we define
\[
A_r = \{X \in SU(n)/SO(n) \mid \lambda_r \text{ is not an eigenvalue of } X\}.
\]
Observe here that we regard $X \in SU(n)/SO(n)$ as a matrix in $SU(n)$ by
 \fullref{homeo}.
Then the $A_r$'s are clearly open sets, and form a covering of $SU(n)/SO(n)$,
 since the property $\lambda_1\lambda_2\cdots\lambda_n \neq 1$ implies that
 the $\lambda_r$'s cannot all appear as the eigenvalues of any matrix in
 $SU(n)/SO(n)$.

Now we fix $A_r$ and let $B$ be a connected component of $A_r$.
In order to show that $A_r$ is contractible in $SU(n)/SO(n)$, it is sufficient
 to show that $B$ is so, since $SU(n)/SO(n)$ is pathwise connected.

Next, let $\u(n) = \{X \in M(n,\C) \mid X^* = -X\}$, and we will define a 
 map $\log \co B \to \u(n)$ as follows.
Let $X \in B \subset A_r$ and $\lambda_r = e^{i\alpha}$, where
 $0 \leq \alpha < 2\pi$.
Then $X$ can be diagonalized by a suitable matrix $P \in U(n)$ as
 $X = PD(e^{i\theta_1},\dots,e^{i\theta_n})P^*$, where $D(a_1,\dots,a_n)$
 denotes a diagonal matrix defined by
\[
D(a_1,\dots,a_n) =
\left(
\begin{array}{@{}c@{}c@{}c@{}}
a_1 &        &     \\[-1.5ex]
    & \ddots &     \\[-1.5ex]
    &        & a_n
\end{array}
\right),
\]
 and we may take $\alpha < \theta_j < \alpha + 2\pi$ for each $j$, since $X$
 does not have $\lambda_r = e^{i\alpha}$ as its eigenvalue.
We define a function $\log\co B \to \u(n)$ by
\[
\log X = PD(i\theta_1,\dots,i\theta_n)P^*,
\]
 where it is easy to see that the definition does not depend on the choice of
 $P$, and the function $\log$ is clearly continuous.
Since $X = \exp(\log X)$ by definition, we have
\[
1 = \det X = \det(\exp(\log X)) = \exp(\tr(\log X)).
\]
Since the maps $\tr\co M(n,\C) \to \C$, which is the trace function, and
 $\log\co B \to \u(n)$ are continuous and since $B$ is connected, there exists
 an integer $k$ such that $\tr(\log X) = 2\pi ik$ for all $X \in B$.

Now we define a constant matrix $X_0$ in $SU(n)/SO(n)$ by
\[
X_0 = \exp\left(\frac{2\pi ik}{n}\right) \cdot E_n,
\]
 and we show that $B$ is contractible to $X_0$.
In order to define a contracting homotopy, we use the fact that $\u(n)$ is
 a vector space, which allows us to construct linear homotopies.
We define a homotopy $F\co B \times [0,1] \to SU(n)/SO(n)$ by
\[
F(X,s) = \exp\left((1-s)\log X + s\frac{2\pi ik}{n}E_n\right).
\]
Clearly, the function $F$ is continuous such that
 $F(X,0) = \exp(\log X) = X$ and $F(X,1) = X_0$ for all $X \in B$.
Here we need to check that $F(X,s) \in SU(n)/SO(n)$ for all $X \in B$ and
 $s \in [0,1]$.
Since $\u(n)$ is the Lie algebra of $U(n)$, we have $F(X,s) \in U(n)$.
Hence it is sufficient to show that $\det(F(X,s)) = 1$ and that
 $\t{F(X,s)} = F(X,s)$.
The former equality can be seen as follows:
$$\eqalignbot{
\det(F(X,s))
&= \det\left(\exp\left((1-s)\log X + s\frac{2\pi ik}{n}E_n\right)\right) \cr
&= \exp\left(\tr\left((1-s)\log X + s\frac{2\pi ik}{n}E_n\right)\right) \cr
&= \exp\left((1-s)\tr(\log X) + s\frac{2\pi ik}{n}\tr(E_n)\right) \cr
&= \exp\left(2\pi ik(1-s) + 2\pi iks\right) \cr
&= \exp\left(2\pi ik\right) \cr
&= 1.}$$
The latter equality can be seen as follows:
$$\eqalignbot{
\t{F(X,s)}
&= \t{\exp\left((1-s)\log X + s\frac{2\pi ik}{n}E_n\right)} \cr
&= \exp\left((1-s)\log\t{X} + s\frac{2\pi ik}{n}\,\t{E_n}\right) \cr
&= \exp\left((1-s)\log X + s\frac{2\pi ik}{n}E_n\right) \cr
&= F(X,s).}\eqno{\qed}$$

\section{L-S category of $SU(2n)/Sp(n)$}

In this section, we will prove (2) of \fullref{main}.
The integral cohomology ring of $SU(2n)/Sp(n)$ is given as follows (see for
 example \cite{MT}):
\[
H^*(SU(2n)/Sp(n);\Z) = \Lambda(x_5,x_9,\dots,x_{4n-3}),
\]
 where $\Lambda$ denotes exterior algebra.
Since the cup-length gives a lower bound of the L-S category, we have
\[
n-1 = \cupl_{\Z}(SU(2n)/Sp(n)) \leq \cat(SU(2n)/Sp(n)).
\]
Thus in order to determine $\cat(SU(2n)/Sp(n))$, it is sufficient to show the
 following proposition.

\begin{prop}\label{main2}
$
\cat(SU(2n)/Sp(n)) \leq n-1
$
\end{prop}

In order to prove \fullref{main2}, we need some lemmas.

\begin{lem}\label{dim}
Let $X$ be any matrix in $SU(2n)/Sp(n)$.
If $\lambda$ is an eigenvalue of $X$, then
\[
\dim W_\lambda \geq 2,
\]
 where $W_\lambda \subset \C^{2n}$ denotes the corresponding eigenspace.
\end{lem}
\begin{proof}
There exists an eigenvector $v \neq 0$ in $\C^{2n}$ such that $Xv = \lambda v$.
Since $X$ satisfies $XX^* = E_{2n}$ and $\t{X} = JX\t{J}$, it follows by an
 easy calculation that $X(J\bar{v}) = \lambda(J\bar{v})$.
Consequently we have that if $v$ is an eigenvector of $\lambda$, so is
 $J\bar{v}$.
Hence it is sufficient to prove that $v$ and $J\bar{v}$ are linearly
 independent.
If $av + bJ\bar{v} = 0\ (a,b \in \C)$, we have
 $\bar{a}\bar{v} + \bar{b}Jv = 0$, and by solving the simultaneous equations,
we see $(|a|^2 + |b|^2)v = 0$, which implies $a = b = 0$.
\end{proof}

Let $\lambda_1,\lambda_2,\dots,\lambda_n$ be different complex numbers with
 $|\lambda_r| = 1$ such that \linebreak
 $\lambda_1^2\lambda_2^2\cdots\lambda_n^2 \neq 1$.
For $1 \leq r \leq n$ we define
\[
A_r = \{X \in SU(2n)/Sp(n) \mid \lambda_r \text{ is not an eigenvalue of } X\}.
\]

\begin{lem}\label{cover}
A family $\{A_r\}_{1 \leq r \leq n}$ forms an open covering of
 $SU(2n)/Sp(n)$\textup{:}
\[
SU(2n)/Sp(n) = \bigcup_{r=1}^n A_r.
\]
\end{lem}
\begin{proof}
Let $X \in (SU(2n)/Sp(n)) \setminus \bigcup_{r=1}^n A_r
 = \bigcap_{r=1}^n \{(SU(2n)/Sp(n)) \setminus A_r\}$.
Then $X$ has $\lambda_r$ as its eigenvalue. Furthermore we see
 by \fullref{dim} that the multiplicity of the eigenvalue $\lambda_r$ is $2$
 for each $r$.
Consequently, $X$ can be diagonalized by a suitable matrix $P \in U(2n)$:
\[
X = PD(\lambda_1,\lambda_1,\dots,\lambda_n,\lambda_n)P^*.
\]
Therefore we have that
 $\det X = \lambda_1^2\lambda_2^2\cdots\lambda_n^2 \neq 1$ which
 contradicts the fact $X \in SU(2n)$.
\end{proof}

\proof[Proof of \fullref{main2}]
By \fullref{cover}, it is sufficient to show that each $A_r$ is contractible
 in $SU(2n)/Sp(n)$; but since $SU(2n)/Sp(n)$ is pathwise connected, it is
 sufficient to show that any connected component of $A_r$ is contractible
 in $SU(2n)/Sp(n)$.
Now we fix $A_r$ and let $B$ be a connected component of $A_r$.
We will show that $B$ is contractible in $SU(2n)/Sp(n)$.

In a similar way to that in \fullref{SU(n)/SO(n)}, we can define a continuous
 function $\log\co B \to \u(2n)$ such that $\exp(\tr(\log X)) = 1$ for $X \in B$.
Then, as was seen before, there exists an integer $k$ such that
 $\tr(\log X) = 2\pi ik$ for all $X \in B$.

Now define a constant matrix $X_0$ in $SU(2n)/Sp(n)$ by
\[
X_0 = \exp\left(\frac{\pi ik}{n}\right) \cdot E_{2n}
\]
 and a contracting homotopy $F\co B \times [0,1] \to SU(2n)/Sp(n)$ by
\[
F(X,s) = \exp\left((1-s)\log X + s\frac{\pi ik}{n}E_{2n}\right).
\]
Clearly, $F$ is continuous such that $F(X,0) = \exp(\log X) = X$ and
 $F(X,1) = X_0$ for all $X \in B$.
Here we need to check that $F(X,s) \in SU(2n)/Sp(n)$ for all $X \in B$ and
 $s \in [0,1]$.
Since $\u(2n)$ is the Lie algebra of $U(2n)$, we have $F(X,s) \in U(2n)$.
Hence it is sufficient to show that $\det(F(X,s)) = 1$ and that
 $\t{F(X,s)} = JF(X,s)\t{J}$.
The former equality can be seen as follows:
\begin{align*}
\det(F(X,s))
&= \det\left(\exp\left((1-s)\log X + s\frac{\pi ik}{n}E_{2n}\right)\right) \\
&= \exp\left(\tr\left((1-s)\log X + s\frac{\pi ik}{n}E_{2n}\right)\right) \\
&= \exp\left((1-s)\tr(\log X) + s\frac{\pi ik}{n}\tr(E_{2n})\right) \\
&= \exp\left(2\pi ik(1-s) + 2\pi iks\right) \\
&= \exp\left(2\pi ik\right) \\
&= 1.
\end{align*}
The latter equality can be seen as follows:
$$\eqalignbot{
\t{F(X,s)}
&= \t{\exp\left((1-s)\log X + s\frac{\pi ik}{n}E_{2n}\right)} \cr
&= \exp\left((1-s)\log\t{X} + s\frac{\pi ik}{n}\,\t{E_{2n}}\right) \cr
&= \exp\left((1-s)\log(JX\t{J}) + s\frac{\pi ik}{n}E_{2n}\right) \cr
&= \exp\left((1-s)J(\log X)\t{J} + s\frac{\pi ik}{n}E_{2n}\right) \cr
&= \exp\left(J\left((1-s)\log X + 
                       s\frac{\pi ik}{n}E_{2n}\right)\t{J}\right) \cr
&= J\exp\left((1-s)\log X + s\frac{\pi ik}{n}E_{2n}\right)\t{J} \cr
&= JF(X,s)\t{J}.}\eqno{\qed}$$

\section{Lusternik--Schnirelmann category of the irreducible symmetric
         Riemann space of classical type}

First let us recall a theorem due to Ganea (for a proof see \cite{Ga}):

\begin{prop}\label{prop_dimension}
If $X$ is an $(r-1)$--connected CW--complex for $r \geq 1$, then
\[
\cat(X) \leq \dim(X)/r.
\]
\end{prop}

We show the following

\begin{prop}\label{cat(V)}
If $V$ is a simply connected, complex $d$--manifold which admits a K\"ahler
 metric, then 
\[
\cat(V) = d.
\]
\end{prop}

In fact, following James \cite{Ja}, it is proved as follows; we have
 $\cat(V) \leq d$ by \fullref{prop_dimension}, since $V$ is simply
 connected.
But with any K\"ahler metric, there exists a closed $2$--form
 on $V$ whose $d$th power is the volume element and so cannot be cohomologous
 to zero.
Hence we have $\cat(V) \geq d$, since the cup-length gives a lower bound of
 the L-S category.

According to Helgason \cite[page 518]{He}, the irreducible symmetric Riemann spaces
 of classical type which has a Hermitian structure are known to be of type
\begin{center}
A\,III,\ \ BD\,I\,($q=2$),\ \ BD\,II\,($n=2$),\ \ D\,III,\ \ C\,I.
\end{center}

Now we also recall from Proposition 4.1 of \cite{He} the following

\begin{prop}\label{Hermitian}
The Hermitian structure of a Hermitian symmetric space is K\"ahlerian.
\end{prop}

It follows from this proposition that the spaces of above type have
 K\"ahler metric.
Hence we see the L-S category of these spaces by \fullref{cat(V)},
 since a Hermitian symmetric space is a complex manifold by definition
 (see \cite[page 372]{He}).
Thus, with \fullref{main}, the L-S category of the irreducible
 symmetric Riemann space of classical type, except that of
 type BD\,I\,($q \neq 2)$, is determined as follows:

{\small\renewcommand{\arraystretch}{1.4}

\begin{center}
\begin{tabular}{|@{\:\,}l@{\:\,}|@{\:\,}l@{\:\,}|@{\:}c@{\:}|@{\:\,}c@{\:\,}|@{\:}c@{\:}|}
\hline
 & \hspace{4.5em}$G/K$ & K\"{a}hler & dimension & $\cat(G/K)$ \\
\hline
A\,I   & $SU(n)/SO(n)\ (n>2)$
       & no & $(n-1)(n+2)/2$ & $n-1$  \\
\hline
A\,II  & $SU(2n)/Sp(n)\ (n>1)$
       & no & $(n-1)(2n+1)$ & $n-1$  \\
\hline
A\,III & $\begin{array}{@{}l@{}} U(p+q)/(U(p) \times U(q)) \\
                           (p \geq q \geq 1)
          \end{array}$
       & yes & $2pq$ & $pq$  \\
\hline
BD\,I  & $\begin{array}{@{}l@{}} SO(p+q)/(SO(p) \times SO(q)) \\
                           (p \geq q \geq 2,\,p+q \neq 4)
          \end{array}$
       & $\begin{array}{@{}c@{\,}l@{}} \text{yes} & (q=2) \\
                                       \text{no}   & (q \neq 2)
          \end{array}$
       & $pq$
       & $\begin{array}{c@{\,}l} p        & (q=2) \\
                                 \text{?} & (q \neq 2)
          \end{array}$ \\
\hline
BD\,II & $SO(n+1)/SO(n)\ (n \geq 2)$
       & $\begin{array}{@{\,}c@{\,}l@{\,}} \text{yes} & (n=2) \\
                                           \text{no}  & (n \neq 2)
          \end{array}$
       & $n$ & $1$  \\
\hline
D\,III & $SO(2l)/U(l)\ (l \geq 4)$
       & yes & $l(l-1)$ & $l(l-1)/2$  \\
\hline
C\,I   & $Sp(n)/U(n)\ (n \geq 3)$
       & yes & $n(n+1)$ & $n(n+1)/2$  \\
\hline
C\,II  & $\begin{array}{@{}l@{}} Sp(p+q)/(Sp(p) \times Sp(q)) \\
                           (p \geq q \geq 1)
          \end{array}$
       & no & $4pq$ & $pq$  \\
\hline
\end{tabular}
\end{center}}

As for the remaining cases;

Firstly, the space of type BD\,II, the real Stiefel manifold $SO(n+1)/SO(n)$,
 is homeomorphic to $S^n$, and hence we have $\cat(SO(n+1)/SO(n)) = 1$.

Secondly, it is known that the space of type C\,II, the symplectic Grassmann
 manifold $Sp(p+q)/(Sp(p) \times Sp(q))$, is 3--connected.
Hence by \fullref{prop_dimension}, we obtain an upper bound
 $\cat(Sp(p+q)/(Sp(p) \times Sp(q))) \leq 4pq/4 = pq$.
It is also known that the cohomology ring of the symplectic Grassmann manifold
 $Sp(p+q)/(Sp(p) \times Sp(q))$ is similar to that of the complex Grassmann
 manifold $U(p+q)/(U(p) \times U(q))$ (see for example \cite{MT}), so we have
 that
 $\cupl(Sp(p+q)/(Sp(p) \times Sp(q))) = \cupl(U(p+q)/(U(p) \times U(q)))$,
 which is given by $pq$, since the cup-length of $U(p+q)/(U(p) \times U(q))$
 is equal to the L-S category of it.
Hence we obtain a lower bound $\cat(Sp(p+q)/(Sp(p) \times Sp(q))) \geq pq$.

Concluding remark: the mod 2 cohomology of type BD\,I\,($q \neq 2)$,
 $SO(p+q)/(SO(p) \times SO(q))$, $p \geq q > 2$, is not known yet.

\appendix
\section{Appendix}

\proof[Proof of \fullref{homeo}]

$(1)$\qua Let $K_n = \{ X \in SU(n) \mid \t{X} = X \}$ and define an action of
 $P \in SU(n)$ on $K_n$ by
\[
P\!\cdot\! X = PX\t{P} \ \ (X \in K_n).
\]
We will show that $X \in K_n$ is represented as follows:
\[
X = P\t{P} = PE_n\t{P} \ \ (P \in SU(n)).
\]
Let $X \in K_n$.
Since $X+\bar{X}$ and $i(X-\bar{X})$ are real symmetric matrices which commute
 with each other, they can be diagonalized by a suitable matrix $B \in SO(n)$:
\[
\t{B}(X+\bar{X})B = D(a_1,\dots,a_n),\ \ \t{B}i(X-\bar{X})B = D(b_1,\dots,b_n).
\]
Then we have
\[
\t{B}XB = D((a_1-ib_1)/2,\dots,(a_n-ib_n)/2),
\]
 where $|(a_k-ib_k)/2| = 1$ for $1 \leq k \leq n$, since $\t{B}XB \in SU(n)$.
Now we can take complex numbers $c_1,\dots,c_n$ such that
 ${c_k}^2 = (a_k-ib_k)/2$ and $c_1 \cdots c_n = 1$.
Then we have $\t{B}XB = CC = C\t{C}$, where $C = D(c_1,\dots,c_n) \in SU(n)$.
 By taking $P = BC$, we have
\[
X = BC\t{C}\t{B} = P\t{P}\ \ (P \in SU(n)),
\]
 which implies that the action is transitive.

On the other hand, the isotropy group at $E_n$ is given by
\[
\{ P \in SU(n) \mid P\t{P} = E_n \}
= \{ P \in SU(n) \mid \bar{P} = P \}
= SO(n).
\]
Since $SU(n)$ is compact, we obtain
\[
SU(n)/SO(n) = \{ X \in SU(n) \mid \t{X} = X \}.
\]

$(2)$\qua There is an embedding $c'\co Sp(n) \to SU(2n)$ defined by
\[
c'(X) =
\left(
  \begin{array}{@{}cc@{}}
     A & -\bar{B} \\
     B & \bar{A}
  \end{array}
\right)
\ \ (X = A + jB)
\]
 such that $Sp(n) = \{ X \in SU(2n) \mid XJ\t{X} = J \}$.
 
Let $L_{2n} = \{ X \in SU(2n) \mid \t{X} = -X \}$ and define an action of
 $P \in SU(2n)$ on $L_{2n}$ by
\[
P\!\cdot\! X = PX\t{P} \ \ (X \in L_{2n}).
\]
We will show that $X \in L_{2n}$ is represented as follows:
\[
X = PJ\t{P} \ \ (P \in SU(2n)).
\]
Let $\lambda$ be an eigenvalue of $X \in L_{2n}$.
Here observe that $|\lambda| = 1$, since $X \in SU(2n)$.
There exists an eigenvector $v \in \C^{2n}$ such that
 $Xv = \lambda v$ and $|v| = 1$.
Since $X$ satisfies $XX^* = E_{2n}$ and $\t{X} = -X$, it follows by an
 easy calculation that $X\bar{v} = -\lambda\bar{v}$.
Let $W$ be the $2$--dimensional subspace of $\C^{2n}$ spanned by
 $v$ and $\bar{v}$.
By repeating this procedure to the orthogonal complement $W^{\perp}$ of
 $W$, we can take consequently an orthonormal basis
 $\{v_1,\bar{v}_1,\dots,v_n,\bar{v}_n\}$ in $\C^{2n}$ such that
\[
Xv_k = \lambda_kv_k, \ \ X\bar{v}_k = -\lambda_k\bar{v}_k\ \ (k=1,\dots,n),
\]
 where $\lambda_1,\dots,\lambda_n$ are eigenvalues of $X$.
Put
\[
w_k = \frac{1}{\sqrt{2}}(v_k + \bar{v}_k),\ \ 
{w_k}' = \frac{-i}{\sqrt{2}}(v_k - \bar{v}_k)\ \ (k=1,\dots,n).
\]
Then $\{w_1,{w_1}',\dots,w_n,{w_n}'\}$ forms an orthonormal basis in 
 $\R^{2n}$ such that
\[
Xw_k = i\lambda_k{w_k}', \ \ X{w_k}' = -i\lambda_kw_k\ \ (k=1,\dots,n).
\]
Thus we have the following:
\[
\t{B}XB = \left(
            \begin{array}{@{}c@{\,}c@{}}
              O & -D(i\lambda_1,\dots,i\lambda_n) \\
              D(i\lambda_1,\dots,i\lambda_n) & O
            \end{array}
          \right),
\]
 where $B = (w_1,{w_1}',\dots,w_n,{w_n}') \in O(2n)$.
Observe that we can choose $B$ in $SO(2n)$ by replacing $\lambda_1$ with
 $-\lambda_1$, if necessary.

Now we take complex numbers $c_1,\dots,c_n$ such that
 ${c_k}^2 = i\lambda_k$ for each $k$, and let
 $C = D(c_1,\dots,c_n,c_1,\dots,c_n) \in U(2n)$.
\renewcommand\theequation{\thesection.\arabic{equation}}
Then we have 
\begin{equation}\label{CJC}
\t{B}XB = CJC = CJ\t{C}.
\end{equation}
We can choose $C$ in $SU(2n)$.
In fact, since $\t{B}XB \in SU(2n)$, we have
\[
\det(\t{B}XB) = i^{2n}\lambda_1^2\cdots\lambda_n^2
= (i^n\lambda_1\cdots\lambda_n)^2 = 1,
\]
so $i^n\lambda_1\cdots\lambda_n = \pm 1$.
Hence $\det(C) = c_1^2 \cdots c_n^2 = i^n\lambda_1\cdots\lambda_n =  \pm 1$.
If $\det(C) = -1$, then replacing $C$ with that multiplied by
$
\left(
  \begin{array}{@{}c@{\,}c@{}}
    D(0,1,\dots,1) & D(1,0,\dots,0) \\
    D(1,0,\dots,0) & D(0,1,\dots,1)
  \end{array}
\right)
$,
 we have $\det(C) = 1$.

By taking $P = BC$, we deduce by \eqref{CJC} that
\[
X = BCJ\t{C}\t{B} = PJ\t{P}\ \ (P \in SU(2n)),
\]
 which implies that the action is transitive.

On the other hand, the isotropy group at $J$ is given by
\[
\{ P \in SU(2n) \mid PJ\t{P} = J \} = Sp(n).
\]
Since $SU(2n)$ is compact, we obtain
\[
SU(2n)/Sp(n) = \{ X \in SU(2n) \mid \t{X} = -X \}.
\]
Further multiplying by $J$, we obtain
$$\eqalignbot{
SU(2n)/Sp(n) &= \{ JX \in SU(2n) \mid \t{X} = -X \} \cr
             &= \{ X \in SU(2n) \mid \t{(\t{J}X)} = -\t{J}X \} \cr
             &= \{ X \in SU(2n) \mid \t{X}J = JX \} \cr
             &= \{ X \in SU(2n) \mid \t{X} = JX\t{J} \}.}\eqno{\qed}
$$

\bibliographystyle{gtart}
\bibliography{link}

\end{document}